\documentclass[a4paper,12pt]{amsart}

\usepackage{amsmath}
\usepackage{amssymb}
\usepackage{amsthm}

\newcommand{\CAT}{\operatorname{CAT}}
\newcommand{\N}{\operatorname{\mathbb{N}}}
\newcommand{\R}{\operatorname{\mathbb{R}}}
\newcommand{\Sph}{\operatorname{\mathbb{S}}}
\newcommand{\ulim}{\operatorname{\text{$\omega$-$\lim$}}}

\theoremstyle{plain}
\newtheorem{thm}{Theorem}[section]
\newtheorem{prop}[thm]{Proposition}
\newtheorem{lem}[thm]{Lemma}

\theoremstyle{definition}
\newtheorem{exmp}[thm]{Example}
\newtheorem{rem}[thm]{Remark}

\begin{document}

\title
[Free construction of CAT(1) spaces]
{Free construction of CAT(1) spaces}

\author
[Koichi Nagano]
{Koichi Nagano} 

\email
[]
{nagano@math.tohoku.ac.jp / nagano@math.uni-bonn.de}

\address
[]
{\endgraf Mathematical Institute, Tohoku University
\endgraf
Aoba, Sendai, Miyagi, 980-8578, Japan}
\curraddr
[]
{Mathematisches Institut, Universit\"{a}t Bonn
\endgraf
Beringstra{\ss}e 1, D-53115 Bonn, Germany}

\date{August 27, 2005}

\thanks{
Partially supported 
by the 2004 JSPS Postdoctoral Fellowships for Research Abroad.}


\begin{abstract}
We construct several non-trivial examples of $\CAT(1)$ spaces
by using the idea of free construction.
\end{abstract}

\maketitle

\section{Introduction}

A spherical (resp.~Euclidean) building $G$ of dimension $n$
has plenty of apartments;
more precisely,
for any two points in $G$
there exists an isometric embedding
from the standard unit $n$-sphere 
(resp.~the Euclidean $n$-space) to $G$
whose image contains the given two points.
Spherical or Euclidean buildings satisfy other rigid metric properties
in the context of the geometry of $\CAT(1)$ spaces
(cf.~\cite{kleileeb}).
In \cite{charlytc, alytchak, balslytc},
we see
several metric characterizations of buildings.

We say that a $\CAT(1)$ space is {\em $n$-round}
if any two points are contained in 
an isometrically embedded unit $n$-sphere.
For instance,
every spherical building of dimension $n$ is $n$-round.
Every $\pi_{\R}$-gon in the sense of \cite{cdbennet}
is $1$-round.
Every hemispherex of dimension $n \ge 2$
in \cite{ballbrin}
is $(n-1)$-round 
(Example \ref{exmp: hemispherex}).
In general, every $n$-round space is geodesically complete
and $(n-1)$-round,
and it has diameter $\pi$.

A geodesically complete $\R$-tree with no graph structure
has plenty of lines.
This example suggests that 
only from the property of plenty of apartments
one can not derive the rigid structure of buildings.
Nevertheless,
it seems to be non-trivial to see
whether an $n$-round space of dimension $n$
without building structure exists.
The main purpose of this paper
is to construct an $n$-round space of dimension $n$
that is not a spherical building.

The Reshetnyak gluing theorem
(\cite{reshetnyak}, cf.~Chapter II.11 in \cite{bridhaef})
helps us to obtain various examples of $\CAT(1)$ spaces.
For our purpose, we need another idea.
The purpose is achieved by mixing
the idea of gluing 
with that of free construction.
The idea of free construction was used by Tits \cite{jacqtits}
to obtain generalized polygons with large automorphism group.

For a $\CAT(1)$ space $X$,
we denote
by $\dim X$ the geometric dimension of $X$
in the sense of \cite{bkleiner}.
By using the idea of free construction,
we obtain the following:

\begin{thm}\label{thm: round}
Let $(X_{\lambda})_{\lambda \in \Lambda}$ 
be a family of $\CAT(1)$ spaces of diameter $\le \pi$.
Then for any $n$ there exists an $n$-round space $Y$
such that
\begin{enumerate}
\item
$\dim Y = \max\{ n, \sup_{\lambda \in \Lambda} \dim X_{\lambda} \}$;
\item
all points in $Y$ are non-manifold points;
\item
each $X_{\lambda}$ is isometrically embedded in $Y$.
\end{enumerate}
In particular,
there exists an $n$-round space of dimension $n$
that is not a spherical building.
\end{thm}

By our free construction, we also see the following.
For a fixed $l \ge \pi$,
let $(X_{\lambda})_{\lambda \in \Lambda}$ 
be a family of $\CAT(1)$ spaces of diameter $\le l$.
Then for any $n$ there exists a $\CAT(1)$ space
in which any two points are contained in 
an isometrically embedded standard $n$-sphere of diameter $l$,
and in which each $X_{\lambda}$ is isometrically embedded
(Theorem \ref{thm: family}).

Next we consider the following condition.
We say that
two points in a metric space
are {\em antipodes}
if they have distance $\ge \pi$.
A $\CAT(1)$ space is {\em $n$-rolling}
if it has at least one pair of antipodes
and if any antipodes are contained 
in an isometrically embedded unit $n$-sphere.
Every $n$-round space is $n$-rolling.
The unit $n$-hemisphere is $(n-1)$-rolling.
Every $n$-rolling space is geodesic and $(n-1)$-rolling, 
and it has diameter $\pi$.
If an $n$-rolling space has dimension $n$,
then it is geodesically complete (cf.~\cite{balslytc}).
It was proved in \cite{balslytc} that
if an $n$-rolling space of dimension $n$
contains a relatively compact open set,
then it is a spherical building.

A $\CAT(1)$ space has dimension $1$
if and only if it is locally an $\R$-tree
of injectivity radius $\ge \pi$.
A $\CAT(1)$ space of dimension $1$
is $1$-rolling if and only if it is $1$-round.
In general, we have:

\begin{thm}\label{thm: rolling}
For any $n \ge 2$,
there exists an $n$-rolling space 
of dimension $n$
that is not $n$-round.
\end{thm}

In fact, for any $n \ge 2$,
such $n$-rolling spaces of dimension $n$ can be obtained abundantly
(Proposition \ref{prop: rollingnonround}). 

We remark that the free construction discussed below
also works in the $\CAT(\kappa)$ setting.
The details are omitted in this paper.

We mention the relation between our notion
and the asymptotic geometry of $\CAT(0)$ spaces.
If $H$ is a locally compact, geodesically complete $\CAT(0)$ space
of dimension $n$,
and if its Tits ideal boundary is $(n-1)$-rolling,
then $H$ is a Euclidean building.
This is substantially a corollary of 
the works of \cite{charlytc, alytchak, balslytc}
(cf.~Section 7 in \cite{balslytc}).

\section{Preliminaries}

\subsection{CAT(1) spaces}

A minimizing geodesic in a metric space
means
a length-minimizing curve joining two points.
Any single point is also assumed to be
a minimizing geodesic.
A metric space is said to be {\em geodesic}
if any two points can be joined by a minimizing geodesic.
A locally geodesic space is {\em geodesically complete}
if any minimizing geodesic is contained 
in a bi-infinite locally minimizing geodesic.

We denote by $\Sph^n$
the standard unit $n$-sphere.

A complete metric space $X$ with distance $d$
is said to be $\CAT(1)$
if any two points $x, y \in X$ with $d(x,y) < \pi$
can be joined by a minimizing geodesic
and if any geodesic triangle in $X$ with perimeter $< 2\pi$
is not thicker than 
its comparison triangle in $\Sph^2$.
If a $\CAT(1)$ space $X$ is geodesically complete,
then for any $x \in X$
we have a point $y \in X$ with $d(x,y) \ge \pi$.
We refer to \cite{bridhaef} for the basics of $\CAT(1)$ spaces
and the historical remarks in geometry of
spaces with curvature bounded above
in the sense of Alexandrov.

We first quote the Reshetnyak gluing theorem
in the following form:

\begin{thm}\label{thm: gluing}
{\em (\cite{reshetnyak}, cf.~Chapter II.11 in \cite{bridhaef})}
Let $(X_{\lambda})_{\lambda \in \Lambda}$
be a family of $\CAT(1)$ spaces
with closed subspaces $C_{\lambda} \subset X_{\lambda}$,
and $X$ a $\CAT(1)$ space
with a closed subspace $C \subset X$.
Assume that
\begin{enumerate}
\item
we have an isometry between $C$ and $C_{\lambda}$;
\item
$C$ is a $\CAT(1)$ space.
\end{enumerate}
Let $Y := X \sqcup (\bigsqcup_{\lambda \in \Lambda} X_{\lambda}) / \sim$ 
be the quotient metric space
obtained by gluing $X_{\lambda}$ to $X$ along $C = C_{\lambda}$
for each $\lambda$.
Then $Y$ is a $\CAT(1)$ space.
\end{thm}

Remark that
a similar statement holds in the $\CAT(\kappa)$ setting.
We use this gluing theorem
only for the cases where
$C$ is a minimizing geodesic of length $\le \pi$
(Subsection 3.2)
and where
$C$ are antipodes
(Subsection 3.3).

Let $X$ be a $\CAT(1)$ space.
For $x, y, z \in X$ with $d(x,y), d(x,z) < \pi$,
the angle at $x$ between $y$ and $z$
is denoted by $\angle_x(y,z)$.
For $x \in X$,
the space of directions at $x$ 
equipped with the interior distance 
induced from $\angle_x$
is denoted by $\Sigma_xX$.
This is also a $\CAT(1)$ space
(\cite{nikolaev}).
The {\em geometric dimension $\dim X$ of $X$}
is defined inductively as follows
(\cite{bkleiner}):
If $X$ is discrete, then $\dim X := 0$;
in the other case,
\[
\dim X := 1 + \sup_{x \in X} \dim \Sigma_xX.
\]
Some geometric study of the geometric dimension
can be seen in \cite{bkleiner}.

We say that a point $x \in X$ is 
{\em strongly singular in $X$}
if $\Sigma_xX$ has at least two connected components.
A manifold point means a point
that has a neighborhood homeomorphic to a Euclidean space.
We denote by $S(X)$
the set of all non-manifold points in $X$.
Every manifold point in a $\CAT(1)$ space
has the space of directions consisting of
a single connected component.
If $x$ is strongly singular in $X$,
then $x \in S(X)$.

\subsection{Poles}

Let $X$ be a geodesically complete $\CAT(1)$ space
of dimension $\ge 2$.
A point $p \in X$ with $\dim \Sigma_pX \ge 1$ is called a {\em pole}
if the diameter of $\Sigma_pX$ is greater than $\pi$.
Remark that
if $X$ is $n$-round for $n \ge 2$,
then it has no pole.

\begin{exmp}\label{exmp: pole}
For $n \ge 2$,
the spherical suspension over 
a standard $(n-1)$-sphere of diameter $> \pi$
is a geodesically complete $\CAT(1)$ space
of dimension $n$ and diameter $\pi$.
The suspension points are poles.
\end{exmp}

\begin{exmp}\label{exmp: twospheres}
Let $n \ge 2$,
and let $p, q \in \Sph^n$ be antipodes.
Take some antipodes $\overline{p}, \overline{q}$
in another unit $n$-sphere $\Sph^n$.
Let $X$ be the quotient metric space
$\Sph^n \sqcup \Sph^n / \sim$
made by attaching $\Sph^n$ to $\Sph^n$ 
at $p = \overline{p}$ and $q = \overline{q}$.
This is a geodesically complete $\CAT(1)$ space
of dimension $n$ and diameter $\pi$.
The antipodes $p$ and $q$ in $X$ are poles. 
\end{exmp}

\begin{exmp}\label{exmp: hemispherex}
Let $n \ge 2$.
A spherical simplicial complex 
is called an {\em $n$-hemispherex}
if it is obtained from $\Sph^n$ 
by attaching unit $n$-hemispheres 
along great-hyperspheres
such that
no pair of antipodes in the original $\Sph^n$
belongs to all hyperspheres (\cite{ballbrin}).
Every $n$-hemispherex is $(n-1)$-round, not $n$-rolling,
and has a pole.
\end{exmp}

\subsection{Ultralimits}

We review the notion of ultralimits,
which goes back to \cite{driewilk}
(cf.~Chapter I.5 in \cite{bridhaef}).

We fix a non-principal ultrafilter $\omega$ on $\N$.
Let $(Y_i)_{i \in \N}$ be a sequence of metric spaces
with distances $d_i$ and basepoints $p_i$.
We denote by $Y_{\omega}^0$
the set of all sequences $(y_i)$
such that $y_i \in Y_i$ and
$d_i(p_i,y_i)$ is uniformly bounded.
Let $d_{\omega} \colon Y_{\omega}^0 \times Y_{\omega}^0 \to [0,\infty)$
be the function defined as
$d_{\omega} ((y_i),(z_i)) := \ulim d_i(y_i,z_i)$,
where $\ulim d_i$ is the ultralimit
of the maps $d_i \colon Y_i \times Y_i \to [0,\infty)$.
We define 
$\ulim Y_i$ as the quotient
metric space $(Y_{\omega}^0,d_{\omega}) / d_{\omega}=0$,
called the {\em ultralimit of $(Y_i)$ 
with respect to $\omega$}.
Note that
if the diameters of $Y_i$ are uniformly bounded,
then $\ulim Y_i$ does not depend on the choices of
the basepoints $p_i$.
If each $Y_i$ is complete
(resp.~geodesically complete),
then $\ulim Y_i$ is again complete (resp.~geodesically complete).
If we have a metric space $X$ and
isometric embeddings $f_i \colon X \to Y_i$,
then we have a natural isometric embedding
$f_{\omega} \colon X \to \ulim Y_i$
defined by
$f_{\omega}(x) := (f_i(x))$.
If each $Y_i$ is $\CAT(1)$,
then $\ulim Y_i$ is again $\CAT(1)$;
moreover,
if $\dim Y_i \le n$ for each $i$,
then $\dim (\ulim Y_i) \le n$
(cf.~\cite{alytchak}).

Let $Y$ be a metric space,
and $\omega$ a non-principal ultrafilter on $\N$.
We call the ultralimit $\ulim Y$
for the constant sequence $(Y)$
the {\em ultracompletion of $Y$} with respect to $\omega$.
By definition, $Y$ can be naturally isometrically embedded in $\ulim Y$.
In the sequel,
we use only ultracompletions of metric spaces. 

\section{Construction}

We prove Theorems \ref{thm: round} and \ref{thm: rolling}.

\subsection{Truncated inductive limits}

For metric spaces $X$ and $Y$,
we write $X < Y$ 
if we have an isometric embedding from $X$ to $Y$.
We say that a map $f \colon X \to Y$ 
is a {\em $\pi$-truncated isometric embedding}
if for any $x, y \in X$ with $d(x,y) \le \pi$ 
we have $d(f(x),f(y)) = d(x,y)$.
We write $X \prec Y$
if we have a $\pi$-truncated 
isometric embedding from $X$ to $Y$.

Let $Y_0 \prec Y_1 \prec \cdots \prec Y_{i-1} \prec \cdots$ 
be an inductive system
of metric spaces $Y_{i-1}$ with distances $d_{i-1}$, $i \in \N$,
and $\pi$-truncated isometric embeddings from $Y_{i-1}$ to $Y_i$.
We denote by $Y_{\infty} := \bigcup_{i \in \N} Y_{i-1}$ 
the inductive limit space
with the distance $d$
defined by $d \, | Y_{i-1} \times Y_{i-1} := \min \{ \pi, d_{i-1} \}$.
We call $Y_{\infty}$ the {\em $\pi$-truncated inductive limit}.
If $Y_i$ is geodesically complete for each $i$,
then so is $Y_{\infty}$.
We need the following later on.

\begin{lem}\label{lem: limit}
For a sequence of $\CAT(1)$ spaces $Y_{i-1}$ of $\dim Y_{i-1} \le m$
with $Y_{i-1} \prec Y_i$, $i \in \N$,
let $Y_{\infty} = \bigcup_{i \in \N} Y_{i-1}$ be the 
$\pi$-truncated inductive limit,
and $\ulim Y_{\infty}$ the ultracompletion of $Y_{\infty}$.
Then 
$\ulim Y_{\infty}$ is a $\CAT(1)$ space of $\dim (\ulim Y_{\infty}) \le m$.
Moreover,
if in addition $Y_i$ is geodesically complete for each $i$,
then we have the following:
\begin{enumerate}
\item
if each point in $Y_{\infty}$ is strongly singular in some $Y_i$,
then each point in $\ulim Y_{\infty}$
is again strongly singular;
\item
if any two points in $Y_{\infty}$ are contained in 
an isometrically embedded unit $n$-sphere,
then $\ulim Y_{\infty}$ is $n$-round;
\item
if any antipodes in $Y_{\infty}$ are contained in 
an isometrically embedded unit $n$-sphere,
then $\ulim Y_{\infty}$ is $n$-rolling.
\end{enumerate}
\end{lem}

\begin{proof}
Since each $Y_i$ is $\CAT(1)$,
$\ulim Y_{\infty}$ is again $\CAT(1)$.
The dimension does not increase under the ultralimits procedures
(cf.~\cite{alytchak}).
Hence we have $\dim (\ulim Y_{\infty}) \le m$.
Observing the ultralimits of
isometrically embedded unit $n$-spheres in $\ulim Y_{\infty}$,
we have (2) and (3).
We show (1).
Suppose the contrary, i.e.,
we have a point $(y)$ in $\ulim Y_{\infty}$
such that $\Sigma_{(y)}(\ulim Y_{\infty})$
has a unique connected component,
where $y \in Y_{\infty}$.
Then $y$ is strongly singular in some $Y_i$.
Hence $\Sigma_yY_i$ has at least two components.
Since each $Y_i$ is geodesically complete,
$\Sigma_{(y)}(\ulim Y_{\infty})$
must have at least two components.
This is a contradiction.
\end{proof}

\subsection{Round spaces}

We prove Theorem \ref{thm: round}.
For a $\CAT(1)$ space $X$,
we denote by $\Gamma_X$
the set of all minimizing geodesics in $X$ of length $\le \pi$.

By using the idea of free construction,
we first show the following:

\begin{prop}\label{prop: roundone}
Let $X$ be a $\CAT(1)$ space of diameter $\le \pi$.
Then for any $n$ there exists an $n$-round space $Y$
with $\dim Y = \max\{ n, \dim X \}$, 
$Y = S(Y)$,
and
$X < Y$.
\end{prop}

\begin{proof}
For a fixed $n$, let $m := \max \{ n, \dim X \}$.
We construct a $\CAT(1)$ space $Y$ inductively as follows.
Put $Y_0 := X$.
Let $Y_{i-1}$ be a 
$\CAT(1)$ space of $\dim Y_{i-1} \le m$.
Set 
$\Gamma_{i-1} := \Gamma_{X_{i-1}}$.
For each $\gamma \in \Gamma_{i-1}$,
we take a minimizing geodesic $\overline{\gamma}$
in $\Sph^n = \Sph^n(\gamma)$ of the same length.
We define $Y_i$ as the quotient metric space
$Y_{i-1} 
\sqcup
\bigl( \bigsqcup_{\gamma \in \Gamma_{i-1}} \Sph^n(\gamma) \bigr) 
/ \sim$
made by gluing $\Sph^n(\gamma)$ to $Y_{i-1}$
along $\gamma = \overline{\gamma}$.
For each $i \in \N$,
$Y_i$ is a geodesically complete, geodesic space.
By the Reshetnyak gluing theorem,
$Y_i$ is $\CAT(1)$.
It follows that $\dim Y_i = m$.
Let $f_i \colon Y_{i-1} \to Y_i$ be the inclusion.
This is a $\pi$-truncated isometric embedding,
and hence $Y_{i-1} \prec Y_i$.
For any $x, y \in Y_{i-1}$,
$f_i(x)$ and $f_i(y)$ are contained in an embedded unit $n$-sphere.
For any $z \in Y_{i-1}$,
$f_i(z)$ is strongly singular in $Y_i$.
Let $Y_{\infty} = \bigcup_{i \in \N} Y_{i-1}$
be the $\pi$-truncated inductive limit.
Each point in $Y_{\infty}$ is strongly singular in some $Y_i$,
and any two points are contained in an embedded unit $n$-sphere.
For a given non-principal ultrafilter $\omega$ on $\N$,
let $Y := \ulim Y_{\infty}$
be the ultracompletion of $Y_{\infty}$.
By Lemma \ref{lem: limit},
$Y$ is an $n$-round space with $\dim Y = m$,
$Y = S(Y)$, and $X < Y$.
\end{proof}

We notice that each $Y_i$ has diameter $> \pi$.
The space $Y_{\infty}$
does not depends
on the choices of minimizing geodesics in $\Sph^n$.

\begin{rem}
We recall that
$\CAT(1)$ spaces in this paper are assumed to be complete.
This seems to be natural
since
almost all significant $\CAT(1)$ spaces are complete.

In the proof of Proposition \ref{prop: roundone},
if one does not require the completeness for the definition
of $\CAT(1)$ spaces,
then $Y_{\infty}$ is a space
with the desired property.
To obtain the completeness,
we need to take the ultracompletion $\ulim Y_{\infty}$.
The usual completion $\overline{Y}_{\infty}$ of $Y_{\infty}$ does not work
for our purpose
since we can not see whether
embedding unit $n$-spheres in $Y_{\infty}$
converge to
the desired spheres in $\overline{Y}_{\infty}$.
\end{rem}

Now we prove Theorem \ref{thm: round}.

\begin{proof}
Let $(X_{\lambda})_{\lambda \in \Lambda}$ be a family 
of $\CAT(1)$ spaces of diameter $\le \pi$.
We fix $n$.
For each $X_{\lambda}$,
by Proposition \ref{prop: roundone} 
we have an $n$-round space $X_{\lambda}'$
with
$\dim X_{\lambda}' = \max \{ n, \dim X_{\lambda} \}$,
$X_{\lambda}' = S(X_{\lambda}')$, 
and $X_{\lambda} < X_{\lambda}'$.

Let 
$m := \max \{ n, \sup_{\lambda \in \Lambda} \dim X_{\lambda} \}$.
We construct a space $Y$ inductively 
in the following way.
Put $Y_0 := X_{\lambda_0}'$ for some $\lambda_0$.
Let $Y_{i-1}$ be a geodesically complete,
geodesic $\CAT(1)$ space of $\dim Y_{i-1} \le m$.
Set $\Gamma_{i-1} := \Gamma_{Y_{i-1}}$.
For each $(\gamma,\lambda) \in \Gamma_{i-1} \times \Lambda$,
take a minimizing geodesic $\overline{\gamma}$
in $X_{\lambda}' = X_{\lambda}'(\gamma)$ of the same length.
Let 
\[
\begin{textstyle}
Y_i := Y_{i-1} 
\sqcup
\bigl( \bigsqcup_{(\gamma,\lambda) \in \Gamma_{i-1} \times \Lambda}
X_{\lambda}'(\gamma) \bigr) 
/ \sim
\end{textstyle}
\]
be the quotient metric space
made by gluing $X_{\lambda}'(\gamma)$ to $Y_{i-1}$
along $\gamma = \overline{\gamma}$ for each $\lambda$.
The space
$Y_i$ is a geodesically complete, geodesic $\CAT(1)$ space
of $\dim Y_i = m$. 
Then $Y_{i-1} \prec Y_i$ by the inclusion
$f_i \colon Y_{i-1} \to Y_i$,
and $X_{\lambda}' < Y_i$ for each $\lambda$. 
For any $x, y \in Y_{i-1}$,
$f_i(x)$ and $f_i(y)$ are contained in an embedded unit $n$-sphere.
For any $z \in Y_{i-1}$,
$f_i(z)$ is strongly singular in $Y_i$.
For a given non-principal ultrafilter $\omega$ on $\N$,
let $Y := \ulim Y_{\infty}$ 
be the ultracompletion
of the $\pi$-truncated inductive limit
$Y_{\infty} = \bigcup_{i \in \N} Y_{i-1}$.
By Lemma \ref{lem: limit},
$Y$ is an $n$-round space with 
$\dim Y = m$ and $Y = S(Y)$.
Moreover, $X_{\lambda} < Y$ for each $\lambda$.
Thus we obtain Theorem \ref{thm: round}.
\end{proof}

\subsection{Rolling spaces}

For a $\CAT(1)$ space $X$,
we denote by $A_X$
the subset of $X \times X$ 
consisting of all pairs of antipodes in $X$.

To prove Theorem \ref{thm: rolling},
it suffices to show:

\begin{prop}\label{prop: rollingnonround}
Let $n \ge 2$,
and let $X$ be a geodesically complete, geodesic $\CAT(1)$ space
of dimension $n$ and diameter $\pi$.
Assume that $X$ has a pole.
Then there exists an $n$-rolling space $Y$ 
of dimension $n$
such that
$Y = S(Y)$, $X < Y$, and $Y$ is not $n$-round.
\end{prop}

\begin{proof}
We first construct a space $Y$ as follows.
Put $Y_0 := X$.
Let $Y_{i-1}$ be a geodesically complete, geodesic
$\CAT(1)$ space of dimension $n$,
and set $A_{i-1} := A_{Y_{i-1}}$.
For each $(p,q) \in A_{i-1}$,
take some antipodes
$\overline{p}, \overline{q}$ in $\Sph^n = \Sph^n(p,q)$.
Define
$Y_i$ as 
$Y_{i-1} \sqcup 
(\bigsqcup_{(p,q) \in A_{i-1}} \Sph^n(p,q) / \sim$
made by attaching $\Sph^n(p,q)$ to $Y_{i-1}$
at $p = \overline{p}$ and $q = \overline{q}$.
The space $Y_i$ is a geodesically complete, geodesic 
$\CAT(1)$ space of dimension $n$.
Then $Y_{i-1} \prec Y_i$ by the inclusion $f_i \colon Y_{i-1} \to Y_i$.
For each $(x,y) \in A_{i-1}$,
the points $f_i(x)$ and $f_i(y)$ are contained 
in an embedded unit $n$-sphere.
For each $i$ and for each $z \in Y_i$,
$f_{i+1}(z)$ is strongly singular in $Y_{i+1}$.
Let $Y_{\infty} = \bigcup_{i \in \N} Y_{i-1}$
be the $\pi$-truncated inductive limit.
Each point in $Y_{\infty}$ is strongly singular in some $Y_i$,
and any antipodes in $Y_{\infty}$ are contained 
in an embedded unit $n$-sphere.
For a given non-principal ultrafilter $\omega$ on $\N$,
let $Y := \ulim Y_{\infty}$ be the ultracompletion of $Y_{\infty}$.
By Lemma \ref{lem: limit},
$Y$ is an $n$-rolling space of dimension $n$
with $Y = S(Y)$ and $X < Y$.
Let $f \colon X \to Y$ be an isometric embedding,
and $p \in X$ a pole.
Then $f(p) \in Y$ is again a pole.
Hence $Y$ is not $n$-round
since every $n$-round space has no pole.
This completes the proof.
\end{proof}

The $\CAT(1)$ spaces in Examples \ref{exmp: pole}, \ref{exmp: twospheres},
and \ref{exmp: hemispherex}
satisfy the assumption in Proposition \ref{prop: rollingnonround}.
Hence we have proved Theorem \ref{thm: rolling}.
\qed

\medskip

\subsection{Modification}

For $l \ge \pi$,
we say that a $\CAT(1)$ space is {\em $(n,l)$-round}
if any two points are contained in 
an isometrically embedded standard $n$-sphere $\Sph_l^n$
of diameter $l$.

Let $X$ be a $\CAT(1)$ space.
For $l \ge \pi$,
we denote by $\Gamma_X^l$
the set of all minimizing geodesics in $X$
of length $\le l$.
In the proof of Theorem \ref{thm: round},
by replacing $\Sph^n$ (resp.~ $\Gamma_X$)
with $\Sph_l^n$ (resp.~ $\Gamma_X^l$),
and by considering $l$-truncated isometric embeddings
and $l$-truncated inductive limits
instead of $\pi$-truncated ones,
we obtain:

\begin{thm}\label{thm: family}
For a fixed $l \ge \pi$,
let $(X_{\lambda})_{\lambda \in \Lambda}$ 
be a family of $\CAT(1)$ spaces of diameter $\le l$.
Then for any $n$ there exists an $(n,l)$-round space $Y$
such that
$\dim Y = \max\{ n, \sup_{\lambda \in \Lambda} \dim X_{\lambda} \}$,
$Y = S(Y)$,
and 
$X_{\lambda} < Y$ for each $\lambda \in \Lambda$.
\end{thm}

\medskip

\noindent
{\bf Acknowledgments}.
I am grateful to Alexander Lytchak for his valuable comment,
to Werner Ballmann for his comment on this work,
and to Linus Kramer
for his information on \cite{cdbennet}.
This paper was written during my stay 
at the University of Bonn.
I would like to thank all geometers in Bonn
for their great hospitality.



\begin{thebibliography}{9999}

\bibitem[BB]{ballbrin}
W. Ballmann and M. Brin,
{\em Diameter rigidity of spherical polyhedra},
Duke Math. J. {\bf 122} (1999), 597--609.

\bibitem[BL]{balslytc}
A. Balser and A. Lytchak,
{\em Building-like spaces},
preprint, 2004;
Math. ArXiv., math.MG/0410437.

\bibitem[B]{cdbennet}
C. D. Bennet,
{\em Twin trees and $\lambda_{\Lambda}$-gons},
Trans. Amer. Math. Soc. {\bf 349} (1997),
Number 5, 2069--2084.

\bibitem[BH]{bridhaef}
M. R. Bridson and A. Haefliger,
{\em Metric spaces of non-positive curvature},
Grund. math. Wiss.,
Volume 319,
Springer-Verlag, 1999.

\bibitem[CL]{charlytc}
R. Charney and A. Lytchak,
{\em Metric characterizations of spherical and Euclidean buildings},
Geom. Topol. {\bf 5} (2001), 521--550. 

\bibitem[DW]{driewilk}
L. van den Dries and A. J. Wilkie,
{\em On Gromov's theorem concerning groups of
polynomial growth and elementary logic},
J. Algebra {\bf 89} (1984), 349--374.

\bibitem[K]{bkleiner}
B. Kleiner,
{\em The local structure of length spaces 
with curvature bounded above},
Math. Z. {\bf 231} (1999), 409--456.

\bibitem[KL]{kleileeb}
B. Kleiner and B. Leeb,
{\em Rigidity of quasi-isometries for symmetric spaces
and Euclidean buildings},
Inst. Hautes \'{E}tudes Sci. Publ. Math. {\bf 86} (1997), 115--197 (1998).

\bibitem[L]{alytchak}
A. Lytchak,
{\it Rigidity of sphrical buildings and spherical joins},
preprint, to appear in Geom. Funct. Anal.

\bibitem[N]{nikolaev}
I. G. Nikolaev,
{\em The tangent cone of an Alexandrov space of curvature $\le K$},
manuscripta math. {\bf 86} (1995), 683--689.

\bibitem[R]{reshetnyak}
Yu. G. Reshetnyak,
{\em On the theory of spaces of curvature not greater than $K$},
Mat. Sb. {\bf 52} (1960), 789--798.

\bibitem[T]{jacqtits}
J. Tits,
{\em Endliche Spiegelungsgruppen, 
die als Weylgruppen auftreten},
Invent. Math. {\bf 43}, (1977),
283--295.

\end{thebibliography}
\end{document}